\documentclass[12pt]{article}
\begin{document}
\title{Generalized Twin Prime Formulas}
\author{H. J. Weber\\Department of Physics\\University of Virginia\\
Charlottesville, VA 22904, U.S.A.}
\maketitle
\begin{abstract}
Based on Golomb's arithmetic formulas, Dirichlet series for two classes of 
twin primes are constructed and related to the roots of the Riemann zeta 
function in the critical strip.  
\end{abstract}
\vspace{3ex}
\leftline{MSC: 11N05, 11M06}
\leftline{Keywords: Golomb's identity, product formula for Dirichlet series} 


\section{Introduction}

Golomb used his arithmetic formula~\cite{sg} 
\begin{eqnarray}
2\Lambda(2a-1)\Lambda(2a+1)=\sum_{d|4a^2-1} \mu(d)\log^2 d,~a\geq 2,
\label{gol}
\end{eqnarray}
where $\Lambda(n)$ is the von Mangoldt function and $\mu(n)$ the M\"obius 
function~\cite{mu},~\cite{hw}, as coefficients of a power series which   
naturally converts to a Lambert series. Due to the lack of an Abelian 
theorem for the latter, no further progress toward a solution of the twin 
prime problem was possible. However, the formula may be used to construct 
twin prime Dirichlet series. Such a method is applied here to twin primes 
$p, p'=p+2D$ for odd $D>0$ first and then even $D>0$. Based on the 
corresponding Golomb identity
\begin{eqnarray}
2\Lambda(2a-D)\Lambda(2a+D)=\sum_{d|4a^2-D^2} \mu(d)\log^2 d,~a\geq a_D,
\label{gen1}
\end{eqnarray}  
for the generalized twin prime numbers $p=2a-D, p'=2a+D.$ The numbers $a_D>0$ 
characterize the first twins that are a distance $2D$ apart. While the running 
median $2a$ between $p=2a-D$ and $p'=2a+D$ for odd $D$ is an even number 
independent of $D,$ whereas $a_D$ depends on $D$ and is more irregular. For 
example, $a_D=2$ for $D=1; a_D=4$ for $D=3; a_D=4$ for $D=5; a_D=5$ for $D=7; 
a_D=7$ for $D=9,$ etc.    

For most even $D,$ such as $2, 4, 8, 10, 14, 16, \ldots$ the median $3(2a-1)$ 
between the twin primes $p=3(2a-1)-D$ and $p'=p=3(2a-1)+D$ that are a distance 
$2D$ apart is again a linear function of the running natural number $a.$ 
Golomb's identity for these cases is  
\begin{eqnarray}
2\Lambda(3(2a-1)-D)\Lambda(3(2a-1)+D)=\sum_{d|9(2a-1)^2-D^2}\mu(d)\log^2 d,~
a\geq a_D.
\label{gen2}
\end{eqnarray}
For $6|D$ the median is more irregular. Therefore, these twin prime numbers 
and those not of the form $(2a-D,2a+D); (3(2a-1)-D,3(2a-1)+D)$ will not be 
considered here.  

The strategy will be to decompose the relevant Golomb identity into two factors 
whose generating Dirichlet series, one acting as the prime number sieve and 
the other to implement the constraint $n=4a^2-D^2$ in $\sum_{d|n}\mu(d)
\log^2 d$ for the twins for odd $D$ and $n=9(2a-1)^2-D^2$ for even $D,$ are 
then used in a product formula for Dirichlet series to construct the twin 
prime Dirichlet series. 
 
\section{Twin Prime Generating Dirichlet Series}

{\it Definition~2.1.} The generating Dirichlet series of one factor of 
the term on the rhs of Golomb's identity~(\ref{gol}) is defined as   
\begin{eqnarray}
Z(s)\equiv\sum_{n=2}^{\infty}\frac{1}{n^s}\sum_{d|n}\mu(d)\log^2 d,
~\Re(s)\equiv\sigma>1. 
\label{l1}
\end{eqnarray}
The series converges absolutely for $\sigma>1.$  Multiplying termwise the
Dirichlet series
\begin{eqnarray}
\frac{d^2}{ds^2}\frac{1}{\zeta(s)}=\sum_{n=2}^{\infty}\frac{\mu(n)\log^2 n}
{n^s}\end{eqnarray}
and $\zeta(s),$ which is justified by absolute convergence of both series 
for $\sigma>1,$ gives Eq.~(\ref{l1}). Carrying out the differentiations 
yields   
\begin{eqnarray}
Z(s)=\zeta(s)\frac{d^2}{ds^2}\frac{1}{\zeta(s)}=2\left(\frac{\zeta'(s)}
{\zeta(s)}\right)^2-\frac{\zeta''(s)}{\zeta(s)}=-\frac{d}{ds}\frac{\zeta'}
{\zeta}(s)+\left(\frac{\zeta'(s)}{\zeta(s)}\right)^2.  
\label{l2}
\end{eqnarray}
Among the poles of the meromorphic function $Z(s)$ are the roots $\rho$ of 
the Riemann zeta function in the critical strip $0<\sigma<1,$ which is 
clear from Eq.~(\ref{l2}). The next two lemmas display analytic 
properties of $Z(s)$ that are needed in Section~6.

{\bf Lemma~2.1.} {\it A pole expansion of $Z(s)$ is given by}
\begin{eqnarray}\nonumber
Z(s)&=&\frac{-1}{(s-1)^2}+\sum_{\rho}\frac{1}{(s-\rho)^2}+\frac{1}{2}\frac{d}
{ds}\frac{\Gamma'(\frac{s}{2}+1)}{\Gamma(\frac{s}{2}+1)}\\&+&\left(1+\frac{
\gamma}{2}-\log 2\pi+\frac{1}{s-1}+\frac{1}{2}\frac{\Gamma'(\frac{s}{2}+1)}
{\Gamma(\frac{s}{2}+1)}-\sum_{\rho}\left(\frac{1}{s-\rho}+\frac{1}{\rho}\right)
\right)^2
\label{zpole}
\end{eqnarray}
{\it with $\gamma=0.57721566\cdots$ the Euler-Mascheroni constant and $\rho$
denoting the zeros of $\zeta(s)$ in the critical strip} $0<\sigma<1.$

{\it Proof.} The pole expansions~\cite{mu},\cite{tit} of the meromorphic
functions $\frac{\Gamma'(s)}{\Gamma(s)}, \frac{\zeta'(s)}{\zeta(s)},$
\begin{eqnarray}\nonumber
\frac{\Gamma'(s)}{\Gamma(s)}&=&-\gamma+\sum_{n=1}^{\infty}\frac{s-1}{n(n+s-1)},
\\\frac{\zeta'(s)}{\zeta(s)}&=&\log 2\pi-1-\frac{\gamma}{2}-\frac{1}{s-1}-
\frac{1}{2}\frac{\Gamma'(\frac{s}{2}+1)}{\Gamma(\frac{s}{2}+1)}+\sum_{\rho}
\left(\frac{1}{s-\rho}+\frac{1}{\rho}\right)
\label{zetap}
\end{eqnarray}
in Eq.~(\ref{l2}) lead to the corresponding pole expansion of $Z(s).~\diamond$

Thus, $Z(s)$ has a simple pole at $s=1$ with the residue
\begin{eqnarray}\nonumber
r_Z(1)&=&2[1-\log 2\pi]+\gamma+\frac{\Gamma'}{\Gamma}(\frac{3}{2})
-2\sum_{\rho}\frac{1}{\rho(1-\rho)},\\
&&\gamma+\frac{\Gamma'}{\Gamma}(\frac{3}{2})=\frac{1}{2}\sum_{n=1}^{\infty}
\frac{1}{n(n+\frac{1}{2})},
\label{res1}
\end{eqnarray}
at $s=\rho$ with the residue
\begin{eqnarray}\nonumber
r_Z(\rho)&=&-2[1-\log 2\pi]-\frac{2}{\rho-1}+\frac{2}{\rho}
+\sum_{n=1}^{\infty}\left(\frac{2}{2n+\rho}-\frac{1}{n}\right),\\
&&\sum_{n=1}^{\infty}\left(\frac{2}{2n+\rho}-\frac{1}{n}\right)=
-\gamma-\frac{\Gamma'}{\Gamma}(1+\frac{\rho}{2}),
\label{res2}
\end{eqnarray}
at $s=-2n$ for $n=1, 2,\ldots$ and double poles at $s=\rho,$ all with
coefficients $2,$ and $s=-2n$ for $n=1,2,\ldots.$

{\bf Lemma~2.2.} {\it The functional equation for $Z(s)$ is given by}
\begin{eqnarray}\nonumber
Z(1-s)&+&\frac{d}{ds}\frac{\Gamma'(s)}{\Gamma(s)}-\frac{(\pi/2)^2}{\cos^2 s
\pi/2}=Z(s)+\left(\frac{\Gamma'(s)}{\Gamma(s)}-\log 2\pi-\frac{\pi}{2}
\tan\frac{s\pi}{2}\right)^2\\&+&2\frac{\zeta'(s)}{\zeta(s)}\left[
\frac{\Gamma'(s)}{\Gamma(s)}-\log 2\pi-\frac{\pi}{2}\tan\frac{s\pi}{2}\right].
\label{zftl}
\end{eqnarray}
   
{\it Proof.} Differentiating the functional 
equation~\cite{mu},\cite{tit},\cite{iv} of $\frac{\zeta'(s)}{\zeta(s)},$
\begin{eqnarray}
-\frac{\zeta'(1-s)}{\zeta(1-s)}=-\log 2\pi-\frac{\pi}{2}\tan\frac{s\pi}{2}+
\frac{\Gamma'(s)}{\Gamma(s)}+\frac{\zeta'(s)}{\zeta(s)},
\label{ftl}
\end{eqnarray}
we obtain
\begin{eqnarray}
\frac{\zeta''(1-s)}{\zeta(1-s)}-\left(\frac{\zeta'(1-s)}{\zeta(1-s)}\right)^2=
\frac{d}{ds}\frac{\Gamma'(s)}{\Gamma(s)}+\frac{\zeta''(s)}{\zeta(s)}-
\left(\frac{\zeta'(s)}{\zeta(s)}\right)^2-\frac{(\pi/2)^2}{\cos^2 s\pi/2}
\end{eqnarray}
and thus the corresponding functional relation for $Z(s).~\diamond$

\section{Constraint Generating Dirichlet Series}

{\it Definition~3.1.} The generating Dirichlet series that implements the 
constraint in Golomb's identity is defined as  
\begin{eqnarray}
Q_D(s)=\sum_{a>[D/2]}^{\infty}\frac{1}{(4a^2-D^2)^s},
\label{qdef}
\end{eqnarray}
with $[D/2]$ the integer part of $D/2.$ It converges absolutely for 
$\sigma>1/2.$ 

{\bf Lemma~3.1.}
\begin{eqnarray}
Q_D(s)=4^{-s}\sum_{\nu=0}^{\infty}(-1)^\nu\left(\frac{D}{2}
\right)^{2\nu}\left(\begin{array}{ll}-s\\ \nu\end{array}\right)
[\zeta(2s+2\nu)-\sum_{a\leq [D/2]}a^{-2s-2\nu}],
\label{qd}
\end{eqnarray}
{\it where $[D/2]$ is the integer part of} $D/2.$

{\it Proof.} Using a binomial expansion and interchanging the summations, 
which is justified by absolute convergence for $\sigma>1/2,$ we obtain 
Eq.~(\ref{qd}).~$\diamond$

Therefore, $Q_D(s)$ has a simple pole at $s=1/2$ with residue $1/4$ and 
at $s=\frac{1}{2}(1-2\nu)$ with residue $\frac{1}{4}(-D)^\nu\left(
\begin{array}{ll}\nu-\frac{1}{2}\\ \nu\end{array}\right)$ for 
$\nu=1, 2, \ldots$ and is regular elsewhere.

{\it Definition~3.2.} The subtracted constraint Dirichlet series is defined as 
\begin{eqnarray}\nonumber
q_D(s)&=&Q_D(s)-2^{-2s}[\zeta(2s)-\sum_{a\leq [D/2]}a^{-2s}]\\&=&
\sum_{a>[D/2]}\left(\frac{1}{(4a^2-D^2)^s}-\frac{1}{(4a^2)^s}\right).
\end{eqnarray}

{\bf Corollary~3.1.} $q_D(s)$ {\it can be resummed as the contour integral}
\begin{eqnarray}
q_D(s)=\frac{-1}{2\pi i}\int_{C}\left(\frac{\Gamma'(2+z)}{\Gamma(2+z)}
-\log z\right)[(4z^2-D^2)^{-s}-(4z^2)^{-s}]dz,~\sigma>\frac{1}{2}, 
\end{eqnarray} 
{\it where the contour $C$ runs parallel to the imaginary axis from 
$c-i\infty$ to $c+i\infty$ with $-[D/2]-1<c<-[D/2],$ and $[D/2]$ the 
integral part of $D/2$.} 
 
{\it Proof.} Adapting a variant of the integral representation of the zeta 
function due to Kloosterman~\cite{tit},   
\begin{eqnarray}\nonumber
\zeta(2s)-\sum_{a\leq [D/2]}a^{-2s}&=&\frac{-1}{2\pi i}
\int_{(-[D/2]-1<c<-[D/2]}\left(\frac{\Gamma'(2+z)}{\Gamma(2+z)}
-\log z\right)z^{-2s}dz,\\&~&\sigma>1/2,
\end{eqnarray}
where the contour runs parallel to the imaginary axis through the abscissa 
$c$. The integral representation can also be obtained by folding the contour 
to the left, running from $-\infty$ back to $-\infty$ enclosing the point 
$-[D/2]-1$ in an anticlockwise sense and applying the residue 
theorem.~$\diamond$ 

{\bf Lemma~3.2.} {\it The Dirichlet series $q_D(s)$ is regular for} 
$\sigma>-1/2,$ ${\cal O}(1)$ {\it for} $|t|\to\infty$ {\it and, with its 
general term grouped as $[(2a)^2-D^2]^{-s}-(2a)^{-2s},$ converges absolutely 
for} $\sigma>-1/2.$ 

{\it Proof.} For $a$ large compared to $|t|,$ the convergence is the same 
as for $t=0$, i.e. $s=\sigma,$ which is a well known property of 
Dirichlet series. A binomial expansion of the general term of $q_D(\sigma),$  
\begin{eqnarray}
(2a)^{-2\sigma}\bigg[\left(1-\frac{D^2}{(2a)^2}\right)^{-\sigma}-1\bigg]
={\cal O}(a^{-2\sigma-2}),
\end{eqnarray} 
shows the absolute convergence and regularity for $\sigma>-1/2.~\diamond$  

The convergence of the integral representation of $q_D(s)$ may be improved 
by including more terms of Stirling's asymptotic series~\cite{aw}. 
This observation leads to 
\begin{eqnarray}\nonumber
q_D(s)&=&\frac{-1}{2\pi i}\int_{C}\left(\frac{\Gamma'(2+z)}{\Gamma(2+z)}
-\log(z+1)-\frac{1}{2(z+1)}+\frac{1}{12(z+1)^2}\right)\\&\cdot&
[(4z^2-D^2)^{-s}-(4z^2)^{-s}]dz.
\label{idc}
\end{eqnarray} 
Now we deform the contour to enclose the real axis in the left-hand plane
along two rays from the origin with opening angle $\phi/|t|$ around $\pi$ for 
fixed $\Im(s)=t\neq 0$ and  $-\pi<\phi<\pi,$ i.e. $z=-re^{i\phi/|t|},~r>0$ is 
used to show that the general term is ${\cal O}(1)$ for $t\to \pm\infty.$ 
The conic area bounded by the rays is closed off by a small circular path 
around $z=-[D/2]-1.$ Stirling's series for the logarithmic derivative of the 
Gamma function applies outside the cone area guaranteeing convergence of the 
integral in Eq.~(\ref{idc}). Each term in
\begin{eqnarray}\nonumber
&&[4z^2-D^2]^{-s}-(2z)^{-2s}\\&=&[(2r)^2e^{2i(\pi-\phi/|t|)}-D^2
]^{-\sigma-it}-[(2r)^2e^{2i(\pi-\phi/|t|)}]^{-\sigma-it}
\end{eqnarray}
remains bounded, in absolute value, as we let $t\to\pm\infty$ for 
$\sigma>-1/2.$ 

\section{Twin Prime Dirichlet Series}

The following product formula for absolutely converging Dirichlet series 
is one of our main tools; it is obtained from summing a well-known formula 
from which mean values of Dirichlet series are usually derived. We mention 
it for ease of reference and because details of its proof are needed in 
Sect.~6. 

{\bf Lemma~4.1.} {\it (Product formula.) Let the Dirichlet series}
\begin{eqnarray}
f(s)=\sum_{m=1}^{\infty}\frac{a_m}{m^s},~g(s)=\sum_{n=1}^{\infty}\frac{b_n}
{n^s}\end{eqnarray}
{\it be single-valued, regular and absolutely convergent for $\sigma>\sigma_a$ 
and $\sigma>\sigma_b,$ respectively. Then the product series}  
\begin{eqnarray}
P(w)=\sum_{l=1}^{\infty}\frac{a_l b_l}{l^w}=\lim_{T\to \infty}\frac{1}{2T}
\int_{-T}^T f(\sigma+it)g(w-\sigma-it)dt
\label{prod}
\end{eqnarray}
{\it is regular and converges absolutely for $\Re(w)\equiv u>\sigma+\sigma_b, 
\sigma>\sigma_a+1.$ The limit of the integral exists and is a regular function 
of the variable} $w\equiv u+iv.$

{\it Proof.} Substituting the Dirichlet series and using~\cite{mu},\cite{tit}
\begin{eqnarray}
\int_{-T}^{T}\left(\frac{n}{m}\right)^{it}dt=\left\{ \begin{array}{ll}2T,~n=m
\\\frac{2\sin(T\log n/m)}{\log n/m},~n\neq m\end{array}\right.,~
\sin\left(T\ln\frac{n}{m}\right)={\cal O}(1),~T\to\infty,
\label{prod1}
\end{eqnarray}
we obtain the first term on the rhs of
\begin{eqnarray}\nonumber
&&\sum_{m,n=1}^{\infty}\frac{a_m b_n}{n^w}\left(\frac{n}{m}\right)^\sigma
\frac{1}{2T}\int_{-T}^T \left(\frac{n}{m}\right)^{it}dt\\&=&
\sum_{n=1}^{\infty}\frac{a_n b_n}{n^w}+\lim_{T\to \infty}
\sum_{1\leq m\neq n}\frac{a_m b_n}{n^w}\left(\frac{n}{m}
\right)^\sigma\frac{\sin(T\log n/m)}{T\log n/m}
=\sum_{n=1}^{\infty}\frac{a_n b_n}{n^w},
\end{eqnarray}
which converges absolutely because $u\geq \sigma_a+\sigma_b+\epsilon.$

To deal with $\sum_{m\neq n}$ we split the summation $m\neq n$ into the
ranges $n\leq m/2, m/2<n<2m, n\geq 2m.$ Since $|\log\frac{n}{m}|\geq 1/\log 2$
for $n\leq m/2$ we obtain the estimate
\begin{eqnarray}
\bigg|\sum_{n\leq m/2}\frac{a_m b_n}{n^u}\left(\frac{n}{m}\right)^\sigma
\frac{\sin(T\log n/m)}{\log n/m}\bigg|\leq\frac{1}{\log 2}\sum_{n\leq m/2}
\frac{|a_m b_n|}{n^{u-\sigma}m^\sigma}={\cal O}(1),
\end{eqnarray}
provided $u-\sigma>\sigma_b$ and $\sigma>\sigma_a.$ For $n\geq 2m,$ we get the
same estimate. For $m/2<n<2m,$ we split the range $\frac{m}{2}<n<2m$ into
$\frac{m}{2}<n<m$ and $m<n<2m$ with $n\neq m$ and use
\begin{eqnarray}
\frac{1}{\log(1-\frac{1}{n})}=n+{\cal O}(1),~\frac{1}{\log(1-\frac{2}{n})}=
\frac{n}{2}+{\cal O}(1),\cdots,~n\sum_{j=1}^{n-1}\frac{1}{j}={\cal O}(n\log n)
\end{eqnarray}
to obtain the estimate
\begin{eqnarray}
\sum_{\frac{m}{2}<n<m}\frac{1}{|\log\frac{n}{m}|}={\cal O}(m\log m),
\end{eqnarray}
and the same estimate for the range $m<n<2m.$ Putting all this together, we
find for $\sigma>\sigma_a+1, u-\sigma>\sigma_b$ that
\begin{eqnarray}\nonumber
\bigg|\sum_{m/2<n<2m}\frac{a_m b_n}{n^u}\left(\frac{n}{m}\right)^\sigma
\frac{\sin(T\log n/m)}{\log n/m}\bigg|&=&{\cal O}\left(\sum_{m/2<n<2m}
\frac{|a_m b_n|\log m}{n^{u-\sigma}m^{\sigma-1}}\right)\\&=&{\cal O}(1).
\end{eqnarray}
Using termwise differentiation with respect to the variable $w$ in the
integrals in conjunction with the estimate $\log n={\cal O}(n^\epsilon)$
shows the absolute convergence of the termwise differentiated product
series and the regularity of the product series.~$\diamond$

As a first step toward constructing the twin prime Dirichlet series we 
apply the product formula to $Z(s)q_D(s-w).$ 

{\bf Theorem~4.1.} {\it For odd} $D>0, \Re(w)=u > \sigma+\frac{3}{2}, 
\sigma>1,$  
\begin{eqnarray}\nonumber
&&\sum_{a>[D/2]}^{\infty}\frac{2\Lambda(2a-D)\Lambda(2a+D)}
{(4a^2-D^2)^w}=A(w)\\&+&\lim_{T\to \infty}\frac{1}{2T}\int_{-T}^T 
Z(\sigma+it)q_D(w-\sigma-it)dt,\label{id4}
\end{eqnarray}
\label{id6}
\begin{eqnarray}
A(w)&\equiv&\sum_{n>[D/2]}^{\infty}\frac{\sum_{d|(2n)^2}\mu(d)
\log^2 d}{(2n)^{2w}},
\label{id5}
\end{eqnarray}
{\it where the limit of the integral on the rhs is a regular function for 
$u>5/2$ and the asymptotic series $A(w)$ converges for $u>1/2.$}

{\it Proof.} Since $Z(s)$ converges absolutely for $\sigma>1,$ the product\\ 
$Z(s)2^{-2w-2s}[\zeta(2w-2s)-\sum_{a\leq [D/2]}a^{-2w+2s}]$ converges 
absolutely for $u-\sigma>1/2,$ i.e. $u>3/2.$ Lemma~4.1 implies that the 
integral in Eq.~(\ref{id6}) is ${\cal O}(T)$ and the limit $T\to\infty$ 
exists representing a regular function for $\sigma>5/2.$ The integral 
over $Z(s)2^{-2(w-s)}[\zeta(2(w-s))-\sum_{a\leq [D/2]}a^{-2w+2s}]$ 
in Eq.~(\ref{id6}) is ${\cal O}(T)$. We apply the product formula 
and obtain the first term on the rhs of Eq.~(\ref{id6}) for 
$u>5/2$.~$\diamond$ 
  
\section{The prime pair generating and asymptotic Dirichlet series}

In order to analyze the generating Dirichlet series $Z(s)$ introduced in 
Section~2 and find the Dirichlet series on the rhs of Eqs.~(\ref{id4}),
(\ref{id5}) in Theor.~4.1, which we call {\it asymptotic Dirichlet series}, 
we evaluate the arithmetic functions in their numerators. 

{\bf Proposition~5.1.} {\it Let $n=\prod_{i=1}^k p_i^{\nu_i}, 
\nu_i\geq 1, P=\prod_{i=1}^k p_i$ be prime number decompositions. Then}
\begin{eqnarray}
\sum_{d|n}\mu(d)\log^2 d=\sum_{d|P}\mu(d)\log^2 d.
\label{p1}
\end{eqnarray}
{\it If $n=p$ is prime then}
\begin{eqnarray}
\sum_{d|p}\mu(d)\log^2 d=-\log^2 p.
\label{p2}
\end{eqnarray}
{\it If $n=p_ip_j$ for prime numbers $p_i\neq p_j$ then}
\begin{eqnarray}
\sum_{d|p_ip_j}\mu(d)\log^2 d=2\log p_i\log p_j.
\label{p3}
\end{eqnarray}
{\it If $k\geq 3$ in $P=\prod_{i=1}^k p_i$ for different prime numbers 
$p_i$ then}
\begin{eqnarray}
\sum_{d|P}\mu(d)\log^2 d=0.
\label{p4}
\end{eqnarray}
{\it Proof.} Eq.~(\ref{p1}) is obvious from the properties of the M\"obius 
function, as is Eq.~(\ref{p2}). Eq.~(\ref{p3}) follows from
\begin{eqnarray}
\sum_{d|p_i p_j}\mu(d)\log^2 d=-\log^2 p_i-\log^2 p_j+\log^2 p_i p_j
=2\log p_i\log p_j.
\end{eqnarray}
Eq.~(\ref{p4}) for $k=3$ follows from expanding $\log^2 p_ip_j=(\log p_i+\log
p_j),^2$ etc.
\begin{eqnarray}\nonumber
&&\sum_{d|p_i  p_j p_l}\mu(d)\log^2 d=-\log^2 p_i-\log^2 p_j-\log^2 p_l
+\log^2 p_i p_j+\log^2 p_i p_l\\\nonumber&+&\log^2 p_j p_l-\log^2 p_i p_j p_l
=2\log p_i\log p_j+2\log p_i\log p_l+\log^2 p_i+\log^2 p_j\\\nonumber&+&
\log^2 p_l+2\log p_j\log p_l-\log^2 p_i-\log^2 p_j-\log^2 p_l\\&-&2\log p_i
\log p_j-2\log p_i\log p_l-2\log p_j\log p_l=0.
\end{eqnarray}
For the general case $k,$ Eq.~(\ref{p4}) is proved by induction. Assuming its
validity for $k$ we can show that for $k+1$:
\begin{eqnarray}\nonumber
0&=&-\log^2 p_{k+1}+\log^2 p_1p_{k+1}+\cdots +\log^2 p_k p_{k+1}\\\nonumber&-&
\log^2 p_1p_2p_{k+1}-\cdots -\log^2 p_{k-1}p_kp_{k+1}+\log^2 p_1p_2p_3p_{k+1}
+\cdots\\&+&\log^2 p_{k-2}p_{k-1}p_kp_{k+1}\pm\cdots +(-1)^{k+1}
\log^2 p_1p_2\ldots p_{k+1}.
\end{eqnarray}
This follows from verifying that the coefficients of $\log^2 p_i$ are zero, as
are those of $2\log p_i p_j.$ In essence, Eqs.~(\ref{p3},\ref{p4}) comprise
Golomb's identity~(\ref{gol}) for $\Lambda(n), \log^2 n.~\diamond$
 
A comment on the constraint $(2a-D,2a+D)=1$ of Golomb's formula is in order.
If $d|2D$ and $d|2a-D$ then $d|2a+D;$ if $d|2a-D,$ but $d$ is not a divisor
of $2D,$ then $d$ is not a divisor of $2a+D.$ Thus, for odd $D,$ if $d|2a-D$
then $d$ does not divide $2a+D,$ except for a {\bf finite} number of divisors
of $D.$ Likewise, if $d|3(2a-1)-D$ and $d$ does not divide $2D$,
then $d$ does not divide $3(2a-1)+D.$ If $a=b\delta$ with $\delta|D$ then,
for odd $D,$ $2a-D=\delta(2b-\frac{D}{\delta}), 2a+D=\delta(2b+\frac{D}
{\delta})$ and $\sum_{d|\delta^2(4b^2-D^2/\delta^2)}\mu(d)\log^2 d=0$
follows from Lemma~5.1. For even $D,$ if $2a-1=b\delta, \delta|D$ with
$\delta\neq 2$ then $9(2a-1)^2-D^2=\delta^2(9b^2-\frac{D^2}{\delta^2})$
and $\sum_{d|\delta^2(9b^2-D^2/\delta^2)}\mu(d)\log^2 d=0$ follows again
from Lemma~5.1. Thus, despite $(2a-D,2a+D)\neq 1, (3(2a-1)-D,3(2a-1)+D)
\neq 1,$ Golomb's identity is trivially satisfied for such values of $a.$

{\bf Lemma~5.1.} {\it The coefficient of the asymptotic Dirichlet series 
$A(w)$ of Theorem~4.1 is given by}
\begin{eqnarray}
\sum_{d|4a^2}\mu(d)\log^2 d=\bigg\{\begin{array}{ll}2\log 2\log p,~a=2^i p^j, 
i\geq 0, j\geq 1\\-\log^2 2, a=2^i, i\geq 0\\0,~a=2^j\prod_{i=1}^k p_i^{\nu_i},
~p_i\neq 2,~k\geq 2,\\\end{array}
\label{mumu}
\end{eqnarray}
{\it where $p, p_i$ are prime numbers $\neq 2$.}

{\it Proof.} This follows from Prop.~5.1. Using the prime number 
decomposition of $a=\prod_{j=1}^k p_j^{\nu_j}, P\equiv \prod_{j=1}^k p_j$ 
in conjunction with 
\begin{eqnarray}
\sum_{d|4a^2}\mu(d)\log^2 d=\sum_{d|2P}\mu(d)\log^2 d,
\end{eqnarray}
the first and third lines of Eq.~(\ref{mumu}) are immediate consequences 
of Golomb's identity. The second line is a special case of $\sum_{d|p}
\mu(d)\log^2 d=-\log^2 p,$ where $p$ is a prime number.$~\diamond$ 

{\bf Theorem~5.1.} {\it For odd $D>0, A(w)$ of Eq.~(\ref{id5}) in 
Theor.~4.1 becomes} 
\begin{eqnarray}\nonumber
A(w)&=&-2\log 2\bigg[(2^{2w}-1)^{-1}\left(\frac{\zeta'}{\zeta}(2w)
+\frac{\log 2}{2^{2w}-1}\right)\\&+&\sum_{2^j p^k\leq [D/2],p>2}
\frac{\log p}{2^{2jw}p^{2kw}}\bigg]-\log^22\bigg[\frac{1}
{2^{2w}-1}-\sum_{1\leq 2^j\leq [D/2]}\frac{1}{2^{2jw}}\bigg].
\label{main}
\end{eqnarray}
{\it Thus, $A(w)$ has a simple pole at $w=1/2$ with the positive 
residue} $\log 2.$ 

{\it Proof.} Summing the Dirichlet series $A(w)$ using Lemma~5.1 
yields 
\begin{eqnarray}\nonumber
&&A(w)=-\log^22\sum_{2^j>[D/2]}\frac{1}{2^{(j+1)2w}}+2\log 2
\sum_{2^j p^k>[D/2],p>2}\frac{\log p}{2^{2jw}p^{2kw}}\\\nonumber
&=&\log^22\bigg[-\frac{1}{2^{2w}-1}+\sum_{1\leq 2^j\leq [D/2]}
\frac{1}{2^{2jw}}\bigg]+2\log 2\sum_{2^j p^k>[D/2],p>2}\frac{\log p}
{2^{2jw}p^{2kw}},\\
\end{eqnarray}
from which Eq.~(\ref{main}) and the residue follow.~$\diamond$ 

\section{Limit of the Integral}

In order to study the integral of Theor.~4.1 we now evaluate the contour 
integral
\begin{eqnarray}
{\cal I}(w)\equiv \frac{1}{2\pi i}\int_C Z(s)q_D(w-s)ds=I_1+I_2+I_3+I_4,
\label{mg}
\end{eqnarray}
where the contour $C$ is a rectangle with the vertices $1+\epsilon-iT, 1+
\epsilon+iT, -\epsilon+iT, -\epsilon-iT$ and $T\neq \gamma, \rho=\beta+i\gamma$ 
the general root of $\zeta(s)$ in the critical strip, $I_1$ is the integral
from $1+\epsilon-iT$ to $1+\epsilon+iT,$ $I_2$ from $1+\epsilon+iT,$ to
$-\epsilon+iT,$ $I_3$ from $-\epsilon+iT,$ to $-\epsilon-iT$ and $I_4$ from
$-\epsilon-iT$ to $1+\epsilon-iT, \epsilon>0.$ The limit $\lim_{T\to\infty}
\frac{\pi}{T}I_1$ of the integral $I_1$ in Eq.~(\ref{id4}) of Theorem~4.1
needs to be analyzed.

{\bf Lemma~6.1.} {\it For} $u>5/2$ 
\begin{eqnarray}\nonumber
{\cal I}(w)&=& q_D(w-1)\bigg[2(1-\log 2\pi)+\sum_{n=1}^{\infty}\frac{1}
{2n(n+\frac{1}{2})}-2\sum_{\rho,|\gamma|<T}\frac{1}{\rho(1-\rho)}\bigg]\\
\nonumber &+&\sum_{\rho,|\gamma|<T}\bigg\{q_D(w-\rho)\bigg[2(\log 2\pi-1)
-\frac{2}{\rho-1}+\frac{2}{\rho}\\&+&\sum_{n=1}^{\infty}\left(\frac{2}
{2n+\rho}-\frac{1}{n}\right)\bigg]-2q'_D(w-\rho)\bigg\}.
\label{rie}
\end{eqnarray}
{\it Proof.} Eq.~(\ref{rie}) follows from using Eq.~(\ref{zpole}) and
the residue theorem taking into account the poles at $s=1, \rho$ 
inside the rectangular contour.

For the evaluation of $I_2, I_4$ we need bounds for $Z(s), q_D(w-s)$ 
for $s=\sigma\pm iT.$

{\bf Lemma~6.2.} $I_2={\cal O}(\log^3 T), I_4={\cal O}(\log^3 T)$ {\it 
for $u>5/2$ and all sufficiently large and and appropriately chosen} $T.$

{\it Proof.} We know~\cite{iv} that for sufficiently large $|t|$ and
$-1<\sigma<2$
\begin{eqnarray}
-\frac{\zeta'(s)}{\zeta(s)}=-\sum_{\rho,|t-\gamma|<1}\left(\frac{1}{s-\rho}+
\frac{1}{\rho}\right)+{\cal O}(\log |t|).
\end{eqnarray}

Choosing $t=\Im(s), T=|t|$ sufficiently large and appropriately we can arrange
that $|t-\gamma|\geq\frac{1}{\log T},$ when $|\gamma-t|<1.$ There are at most
${\cal O}(\log T)$ roots $\rho$ with $|\gamma-t|<1.$ As a result\cite{iv},
\begin{eqnarray}
\sum_{\rho}\frac{1}{|s-\rho|^2}=\sum_{\rho,|\gamma-T|\geq 1}\frac{1}
{|s-\rho|^2}+\sum_{\rho,|\gamma-T|<1}\frac{1}{|s-\rho|^2}={\cal O}(\log T)+
{\cal O}(\log^3 T).
\end{eqnarray}
Hence $|Z(s)|={\cal O}(\log^3 T)$ for $I_2$ and the same bound holds for
$I_4$ and $-\epsilon\leq\sigma\leq 1+\epsilon.$

On all four sides of the rectangular path $-\epsilon\leq \sigma\leq 1+
\epsilon,~\Re(w-s)>-\frac{3}{2}$ for $u>3/2$ the Dirichlet series 
$q_D(w-s)$ remains bounded for $u-\sigma>1/2,$ in absolute magnitude, for 
$|t|\to\infty$ according to Lemma~3.2 .~$\diamond$

{\bf Lemma~6.3.} {\it For} $u>5/2,~I_3={\cal O}(\log^2 T).$

{\it Proof.} For $I_3$ we will use the functional equation~(\ref{zftl}) of
$Z(s)$ to evaluate the first term on the rhs of
\begin{eqnarray}\nonumber
I_3&=&-\frac{1}{2\pi i}\int_{-\epsilon-iT}^{-\epsilon+iT}q_D(w-s)\bigg\{
Z(1-s)+\frac{d}{ds}\frac{\Gamma'(s)}{\Gamma(s)}-\frac{(\pi/2)^2}{
\cos^2 s\pi/2}\\\nonumber&-&\left(\frac{\Gamma'(s)}{\Gamma(s)}-\log 2\pi
-\frac{\pi}{2}\tan\frac{s\pi}{2}\right)^2\\&-&2\frac{\zeta'(s)}{\zeta(s)}
\left[\frac{\Gamma'(s)}{\Gamma(s)}-\log 2\pi-\frac{\pi}{2}\tan\frac{s\pi}{2}
\right]\bigg\}ds,
\label{i3}
\end{eqnarray}
$q_D(w-s)Z(1-s),$ by substituting the absolutely converging Dirichlet series.
This yields a double series $\sum_{m,n}$ as for $I_1$ in Lemma~4.1 involving
the integral $\int_{-T}^{T}(m[(2n)^2-D^2)]^{it}dt.$ Thus, there is no 
contribution for $4n^2-D^2=m$ that would be proportional to $T.$ As earlier, 
$\sum_{m\neq (2n)^2-D^2}$ involves $\sin(T\log [(2n)^2-D^2]m)$ and is 
${\cal O}(1)$ for $u+\epsilon>3/2.$ The same applies to the double sums 
involving $(2n)^2$ and $4n^2-2D^2$ instead of $(2n)^2-D^2.$

To deal with the term involving $-\frac{\zeta'(s)}{\zeta(s)}2(\gamma+\log 2\pi)
q_D(w-s)$ we use the functional equation of Lemma~2.2 and then substitute the
absolutely converging Dirichlet series for $q_D(w-s)$ and $-\frac{\zeta'(1-s)}
{\zeta(1-s)}$ at $\sigma=-\epsilon.$ Again, each double series is obtained as
for $I_1$ in Lemma~4.1 involving the integral $\int_{-T}^{T}(m[(2n)^2-D^2])^{
it}dt,$ $\int_{-T}^{T}(m(2n)^2)^{it}dt,$ or $\int_{-T}^{T}[m(4n^2-2D^2)]^{it}
dt,$ respectively. Again, there is no contribution for $(2n)^2-D^2=m,$
$(2n)^2-2D^2=m,$ or $(2n)^2=m,$ that would be proportional to $T.$ As earlier,
$\sum_{m\neq (2n)^2-D^2},$ $\sum_{m\neq (2n)^2-2D^2},$ or $\sum_{m\neq (2n)^2}$
 lead to the bound ${\cal O}(1).$

Putting all this together, we obtain
\begin{eqnarray}\nonumber
I_3&=&-\frac{1}{2\pi i}\int_{-\epsilon-iT}^{-\epsilon+iT}q_D(w-s)\\\nonumber
&\cdot&\bigg\{\frac{d}{ds}\frac{\Gamma'(s)}{\Gamma(s)}-\frac{(\pi/2)^2}
{\cos^2 s\pi/2}-\left(\frac{\Gamma'(s)}{\Gamma(s)}-\log 2\pi-\frac{\pi}{2}
\tan\frac{s\pi}{2}\right)^2\\\nonumber &+&\left[-2\log 2\pi-\pi\tan\frac{s\pi}
{2}+2\frac{\Gamma'(s)}{\Gamma(s)}\right]\left[\frac{\Gamma'(s)}{\Gamma(s)}
-\log 2\pi-\frac{\pi}{2}\tan\frac{s\pi}{2}\right]\\&+&2\frac{\zeta'(1-s)}
{\zeta(1-s)}\left[\frac{\Gamma'(s)}{\Gamma(s)}+\gamma-\frac{\pi}{2}
\tan\frac{s\pi}{2}\right]\bigg\}ds+{\cal O}(1).
\label{hilf}
\end{eqnarray}
The term $\frac{1}{2\pi i}\int_{-\epsilon-iT}^{-\epsilon+iT}q_D(w-s)
\frac{(\pi/2)^2}{\cos^2 s\pi/2}ds={\cal O}(1)$ because 
$(\cos^2 s\pi/2)^{-2}={\cal O}(e^{-|t|\pi})$ for $t\to \pm\infty.$

The terms involving $|\tan\frac{s\pi}{2}|\to 1$ for $t\to\pm\infty,$ so
$|\frac{d}{ds}\tan\frac{s\pi}{2}|={\cal O}(e^{-|t|\pi/2}).$ Using integration
by parts, we get the following estimate for this term from the Dirichlet
series for $q_D(w-s)$
\begin{eqnarray*}
&&\int_{-\epsilon-iT}^{-\epsilon+iT}\tan\frac{s\pi}{2}\sum_{m=2}^{\infty}
[(4m^2-D^2)^{s-w}-(4m^2)^{s-w}]ds\\&=&\sum_{m=2}^{\infty}[(4m^2-D^2)^{s-w}
(\log(4m^2-D^2))^{-1}\\&-&(4m^2)^{s-w}(2\log(4m^2))^{-1}]
\tan\frac{s\pi}{2}\bigg|_{s=-\epsilon-iT}^{s=-\epsilon+iT}\\&-&
\int_{-\epsilon-iT}^{-\epsilon+iT}\sum_{m=2}^{\infty}[(2m)^2-D^2]^{s-w}
[\log((4m^2-D^2))^{-1}\\&-&(4m^2)^{s-w}
(2\log(4m^2))^{-1}\\&-&(4m^2-2D^2)^{s-w}(2\log(4m^2-2D^2))^{-1}]\frac{d}{ds}
\tan\frac{s\pi}{2}ds,
\end{eqnarray*}
where the integrated term is ${\cal O}(1)$ because the Dirichlet series is
absolutely convergent, $\tan\frac{s\pi}{2}={\cal O}(1)$ and the integral is
${\cal O}(1).$ Terms involving $\tan^2\frac{s\pi}{2},\\\tan\frac{s\pi}{2}
\frac{\Gamma'(s)}{\Gamma(s)}$ are handled in the same way.

The terms involving $q_D(w-s)\frac{\Gamma'(s)}{\Gamma(s)}$ are also treated by
integration by parts where the integrated term is ${\cal O}(\log T)$ and the
integral is integrated by parts leading to $\frac{d}{ds}\frac{\Gamma'(s)}
{\Gamma(s)}={\cal O}(T^{-1})$ at $s=-\epsilon\pm iT.$ Hence the remaining
integral is ${\cal O}(\log T).$ The term involving $q_D(w-s)\frac{d}{ds}
\frac{\Gamma'(s)}{\Gamma(s)}$ is treated by integration by parts leading to
the same estimate. The term $q_D(w-s)\left(\frac{\Gamma'(s)}{\Gamma(s)}
\right)^2$ is also treated by integrating by parts leading to an
${\cal O}(\log^2 T)$ estimate.

The product $q_D(w-s)\frac{\zeta'(1-s)}{\zeta(1-s)}$ is another absolutely
converging Dirichlet series without constant term, i.e. unity. All these
terms in Eq.~(\ref{hilf}) are also treated by integration by parts leading to
the estimate ${\cal O}(\log T).$ Putting all this together proves Lemma~6.3.~
$\diamond$

{\bf Theorem~6.1.} {\it For odd} $D>0, u>5/2,$
\begin{eqnarray}\nonumber
&&\lim_{T\to\infty}\frac{{\pi\cal I}(w)}{T}=\lim_{T\to\infty}\frac{\pi I_1}
{T}=\lim_{T\to\infty}\frac{\pi}{T}\sum_{\rho,|\gamma|<T}\bigg\{q_D(w-\rho)
\cdot\\\nonumber&&\bigg[2(\log 2\pi-1)-\frac{2}{\rho-1}+\frac{2}{\rho}
+\sum_{n=1}^{\infty}\left(\frac{2}{2n+\rho}-\frac{1}{n}\right)\bigg]
-2q'_D(w-\rho)\bigg\},\\&&\lim_{T\to\infty}\frac{\pi I_1}{T}
=\sum_{a>[D/2]}\frac{2\Lambda(2a-D)\Lambda(2a+D)}{(4a^2-D^2)^w}.
\label{59}
\end{eqnarray}

{\it Proof.} We apply the limit $T\to\infty$ on ${\cal I}(w)$ to
Eqs.~(\ref{mg}),(\ref{rie}) using the bounds on $I_2, I_3, I_4$ from
Lemmas~6.2, 6.3 to obtain Eq.~(\ref{59}) because the first line of
Eq.~(\ref{rie}) and the first term of the second line drop out. The 
last line follows from Theor.~4.1 and Lemma~4.1.~$\diamond$

This is our main result. Due to the limit $T\to\infty$ in Theor.~6.1  
the twin prime distributions depend on the asymptotic properties of
the roots of the Riemann zeta function. This feature contrasts with
the analytic proof of the prime number theorem, where the remainder
term is linked to the roots of the Riemann zeta function producing
the staircase-like corrections of the smooth asymptotic limit from
{\bf all} zeta function roots, whereas the leading asymptotic term
has nothing to do with them, originating from the simple pole of
the Riemann zeta function.

\section{Twin Primes For Even $D$}

{\it Definition~7.1.} The corresponding constraint generating 
Dirichlet series is defined as
\begin{eqnarray}
Q_D(s)=\sum_{a>[D/6]}\frac{1}{[3^2(2a-1)^2-D^2]^s},~\sigma>1/2, 
\end{eqnarray}
with $[D/6]$ the integer part of $D/6.$

{\bf Lemma~7.1.} {\it The expansion corresponding to Lemma~3.1 becomes}
\begin{eqnarray}\nonumber
Q_D(s)&=&\sum_{\nu=0}^{\infty}(-1)^\nu\left(\frac{D}{3}\right)^{2\nu}3^{-2s}
\left(\begin{array}{ll}-s\\ \nu\end{array}\right)\bigg[\left(1-\frac{1}
{2^{2s+2\nu}}\right)\zeta(2s+2\nu)\\&-&\sum_{a\leq [D/6]}\frac{1}
{(2a-1)^{2s+2\nu}}\bigg].
\label{qdo}
\end{eqnarray}
{\it Definition~7.2.} The subtracted constraint Dirichlet series is 
defined as
\begin{eqnarray}
q_D(s)=Q_D(s)-3^{-2s}\bigg[(1-2^{-2s})\zeta(2s)-\sum_{a\leq [D/6]}
(2a-1)^{-2s}\bigg].
\end{eqnarray}
Equation~(\ref{idc}) of Lemma~3.2 for this case becomes 
\begin{eqnarray}\nonumber
q(s)&=&\frac{-1}{2\pi i}\int_{C}\left(\frac{\Gamma'(\frac{z-1}{2})}
{\Gamma(\frac{z-1}{2})}-\log(\frac{z-3}{2})-\frac{1}{z-3}+\frac{1}
{3(z-3)^2}\right)\\&\cdot&\bigg[(9z^2-D^2)^{-s}-(3z)^{-2s}\bigg]dz,
\end{eqnarray}
with $-2[D/6]<c<-2[D/6]+1.$ Lemma~3.2 is valid for this case with the 
general term grouped as 
\begin{eqnarray}
\frac{1}{[3^2(2a-1)^2-D^2]^s}-3^{-2s}(2a)^{-2s}.
\end{eqnarray} 
{\bf Theorem~7.1.} {\it For even} $D>0, \Re(w)=u>\sigma+\frac{3}{2}, 
\sigma>1,$  
\begin{eqnarray}\nonumber
&&\sum_{a>[D/6]}^{\infty}\frac{2\Lambda(3(2a-1)-D)\Lambda(3(2a-1)+D)}
{[9(2a-1)^2-D^2]^w}=A(w)
\\&+&\lim_{T\to \infty}\frac{1}{2T}\int_{-T}^T Z(\sigma+it)q(w-\sigma-it)
dt,\label{id4}
\end{eqnarray}
{\it with} 
\begin{eqnarray}
A(w)&\equiv&3^{-2w}\sum_{n>[D/6]}^{\infty}\frac{\sum_{d|3^2(2n-1)^2}
\mu(d)\log^2 d}{(2n-1)^{2w}}.
\end{eqnarray}
This is Theor.~4.1 for even $D.$

{\bf Lemma~7.2.} {\it The general coefficient of the asymptotic 
Dirichlet series $A(w)$ is given by}
\begin{eqnarray}
\sum_{d|9(2a-1)^2}=\bigg\{\begin{array}{ll}-\ln^2 3,~2a-1=3^i\\
2\ln 3\log p,~2a-1=p^i,~p\neq 3\\0,~\rm{~else.~}\\\end{array}
\end{eqnarray}

{\bf Theorem~7.2.} {\it For even $D>0, A(w)$ of Eq.~(\ref{id5}) in 
Theor.~4.1 becomes}
\begin{eqnarray}\nonumber
A(w)&=&-\log^2 3\sum_{3^j>[D/6]}3^{-2jw}+2\log 3\sum_{3^{2jw}p^{2kw}>[D/6], 
p\neq 3}\frac{\log p}{3^{2jw}p^{2kw}}\\\nonumber&=&-\log^23\bigg[
\frac{1}{3^{2w}-1}-\sum_{1\leq 3^j\leq [D/6]}\frac{1}{3^{2jw}}\bigg]\\
\nonumber&-&2\log 3\bigg[(3^{2w}-1)^{-1}\left(\frac{\zeta'}{\zeta}(2w)
+\frac{\log 3}{3^{2w}-1}\right)\\&+&\sum_{3^j p^k\leq [D/6],p\neq 3}
\frac{\log p}{3^{2jw}p^{2kw}}\bigg].
\end{eqnarray} 
Thus, $A(w)$ has a simple pole at $w=1/2$ with the residue $\log 3.$ 

Theor.~6.1 remains valid, except for replacing the twin Dirichlet series in 
Eq.~(\ref{59}) by the lhs of Eq.~(\ref{id4}) in Theor.~7.1. 

\section{Discussion and Conclusion}

The twin prime formulas link the twin primes with the roots of the 
Riemann zeta function, which may be viewed as a step in Riemann's 
program linking prime numbers to these roots. Their dependence on 
the limit $T\to \infty$ demonstrates that the twin prime 
distributions depend on the asymptotic properties of these roots. 
This qualitative feature differs fundamentally from the prime 
number theorem for ordinary prime numbers where all roots 
contribute to the remainder term.

\end{document}